\documentclass[11pt,a4paper]{article}%
\usepackage{amsmath,amsthm}
\usepackage{indentfirst}
\usepackage{cite}
\newtheorem{thm}{Theorem}[section]

\newtheorem{rem}{Remark}[section]

\begin{document}

\title{A simple proof of monotonicity for remainder of Stirling's formula}
\author{Yuling Xue, Songbai Guo\thanks{Corresponding author. E-mail
address: guosongbai@bucea.edu.cn.}}
\date{}
\maketitle
\vskip-6mm
\begin{tabular}[c]{l}
{School of Science, Beijing University of Civil Engineering and
Architecture, Beijing}\\
{\small 102616, P.~R.~China}
\end{tabular}


\begin{quotation}
\noindent\textbf{Abstract} The monotonicity properties of remainder of Stirling's formula for the gamma function are simply obtained by using the integral transforms with series.


\noindent\textbf{AMS(2020)Subject Classification}\quad 33B15, 34A25, 40A25.
\end{quotation}

\section{Introduction}

The famous Stirling's formula for the gamma function is
\begin{equation}
\Gamma(x+1)=\sqrt{2\pi x}\left(  \frac{x}{e}\right)  ^{x}e^{\frac{\sigma(x)}{12x}}, \sigma(x)\in(0,1), x>0, \label{eq1}%
\end{equation}
which is an important formula in fields such as applied mathematics, numerical mathematics and probability.
The related characterization of formula \eqref{eq1} has been drawn much
attention (see, e.g., \cite{Dutkay13,Mortici11,Shi06,Sevli11,Koumandos09,Lin14}). The Binet's first
formula \cite{Magnus66} for the gamma function is as follows:
\begin{equation}
\ln\Gamma(x+1)=\left(  x+\frac{1}{2}\right)  \ln x-x+\frac{1}{2}\ln(2\pi)+\int%
_{0}^{\infty}\phi(t)e^{-xt}dt,\ x>0, \label{eq2}%
\end{equation}
where $$\phi(t)=\left(  \frac{1}{e^{t}-1}-\frac{1}{t}+\frac{1}{2}\right)
\frac{1}{t}.$$ Furthermore, by \cite{Knopp54}, it holds
\begin{equation}
\phi(t)=\sum\limits_{v=1}^{\infty}\frac{2}{4\pi^{2}v^{2}+t^{2}}.
\label{eq3}%
\end{equation}
It follows from formulae \eqref{eq1}-\eqref{eq3} that
\begin{equation}
\sigma(x)=12x\ln\frac{\Gamma(x+1)}{\sqrt{2\pi x}\left(  x/e\right)^{x}}=12x\int_{0}^{\infty}\phi(t)e^{-xt}dt. \label{eq4}%
\end{equation}
Set
\begin{equation*}
\lambda(x)=\frac{\Gamma(x+1)}{\sqrt{2\pi x}\left(  x/e\right)  ^{x}}-1.
\end{equation*}

Recently, the authors in \cite{Guo14,Guo15,Qi17} used different methods to
study on the monotonicity of the function $\sigma(x)$ or $\lambda(x)$, respectively.
Although the monotonicity results on the functions $\sigma(x)$ and $\lambda(x)$ have been obtained,
we find a more direct and simpler proof to get them motivated by \cite{Guo14,Qi17}.
Specifically, by using the integral transforms with series, we prove simply that the function $\sigma(x)$ strictly increases on $(0,\infty)$. Along the line of the proof, we also show directly that $\lambda(x)$ strictly decreases on $(0,\infty)$, and $\theta(x):=1-\sigma(x)$ is strictly completely monotonic on $(0,\infty)$, which also implies the existing related results.

\section{Main Results}
\setcounter{equation}{0}
Now we are ready to state and prove our work.

\begin{thm}\label{thm1}
The function $\sigma(x)$ strictly increases on $(0,\infty)$.
\end{thm}
\proof It follows easily that the series $\sum\limits_{v=1}^{\infty}2/(4\pi^{2}v^{2}+t^{2})$
uniformly converges to $\phi(t)$ on $[0,\infty)$ since
\begin{equation}\label{eq5}
\sum\limits_{v=1}^{\infty}\frac{2}{4\pi^{2}v^{2}+t^{2}}\leq\sum\limits_{v=1}^{\infty}\frac{1}{2\pi^{2}v^{2}}=\frac{1}{12}.
\end{equation}
Thus, $\phi(t)$ is continuous on $[0,\infty)$. By formula \eqref{eq4}, we have
\begin{equation*}
\sigma(x)=12\int_{0}^{\infty}\phi\left(  \frac{s}{x}\right) e^{-s} ds.
\end{equation*}
Formula \eqref{eq3} hints that $\phi(t)$ strictly decreases on $(0,\infty)$.
As a result, for $s>0$, $\phi(s/x)  =\sum\limits_{v=1}%
^{\infty}2/[4\pi^{2}v^{2}+(s/x)^{2}]$ strictly increases
with respect to $x$ on $(0,\infty)$. Whenever $0<x_{1}<x_{2}$, it is not
difficult to verify that
\begin{equation*}
\int_{0}^{\infty}\phi\left(  \frac{s}{x_{1}}\right) e^{-s} ds<\int%
_{0}^{\infty}\phi\left(  \frac{s}{x_{2}}\right) e^{-s} ds. \label{eq7}%
\end{equation*}
Therefore, $\sigma(x)$ strictly increases on $(0,\infty)$.

\begin{thm}
The function $\lambda(x)$ strictly decreases on $(0,\infty)$.
\end{thm}
\proof Denote $h(x)=\ln(\lambda(x)+1)=\sigma(x)/12x$. From
\eqref{eq4}, it follows that $h(x)=\int^{\infty}_{0}\phi(t)e^{-xt}dt$.
It is easy to find that $h(x)$ strictly decreases on $(0,\infty)$ as
$e^{-xt}$ strictly decreases with respect to $x$ on
$(0,\infty)$ for each $t>0$. Thus, $\lambda(x)$ strictly decreases on $(0,\infty)$.

\begin{thm}\label{thm3}
The derivatives of all orders of the function $\theta(x)$ satisfy
\begin{equation}\label{eq8}
(-1)^{n}\theta^{(n)}(x)>0,x>0,n=0,1,2,\cdots.
\end{equation}
\end{thm}
\proof Because
$$\frac{4t}{(4\pi^{2}v^{2}+t^{2})^{2}}=\frac{1}{4\pi^{2}v^{2}+t^{2}}
\frac{4t}{4\pi^{2}v^{2}+t^{2}}<\frac{1}{4\pi^{2}v^{2}+t^{2}}\leq\frac{1}{4\pi^{2}v^{2}},$$
it follows from \eqref{eq5} that the series $\sum\limits_{v=1}^{\infty}4t/(4\pi^{2}v^{2}+t^{2})^{2}$ uniformly converges in $t\geq0$.
By formulae \eqref{eq3}, \eqref{eq4} and the theorem on term-by-term differentiation of series,
it holds
\begin{equation*}
\sigma(x)=-12\int_{0}^{\infty}\phi(t)de^{-xt}=1-12\int_{0}^{\infty}%
\psi(t)e^{-xt}dt,
\end{equation*}
where $\psi(t)=\sum\limits_{v=1}^{\infty}4t/(4\pi^{2}v^{2}+t^{2})^{2}.$
That is, the nonconstant function
\[
\theta(x)=12\int_{0}^{\infty}\psi(t)e^{-xt}dt.
\]
Now use the Hausdorff-Bernstein-Widder theorem \cite{Widder46}.
\begin{rem}
The function $\theta(x)$ equipped with \eqref{eq8} is said to be strictly completely monotonic,
see the definition of Widder \cite[Chap. IV]{Widder46}.
\end{rem}
\begin{rem} Theorem \ref{thm3} implies the main results in \cite{Guo14,Guo15,Qi17}
and the proof of that is also simpler than those in the three papers. It is also not difficult to find that
Theorem \ref{thm1} is a corollary of Theorem \ref{thm3}.
\end{rem}


%

\section*{Acknowledgements}

This work is supported in part by the National Natural Science Foundation of
China (No. 11901027), the China Postdoctoral Science Foundation (No. 2021M703426),
the Pyramid Talent Training Project of BUCEA (JDYC20200327).


\begin{thebibliography}{99}                                                                                               %
\bibitem {Dutkay13}D. E. Dutkay, C. P. Niculescu, F. Popovici, Stirling's formula
and its extension for the gamma function, \textit{Amer. Math. Monthly}
\textbf{120} (8) (2013) 737--740.

\bibitem {Guo14}S. Guo, W. Ma, B. G. S. A. Pradeep, Complete characterizations
of the gamma function, \textit{Appl. Math. Comput.} \textbf{244} (2014) 912--916.

\bibitem {Guo15}S. Guo, Y. Shen, X. Li, A note on complete monotonicity of the
remainder in Stirling's formula, \textit{Math. Notes} \textbf{97} (6) (2015) 961--964.

\bibitem {Knopp54}K. Knopp, \textit{Theory and Application of Infinite Series},
2nd ed., Blackie \& Son Limited, Glasgow, 1954: 378.

\bibitem {Koumandos09}S. Koumandos, H. L. Pedersen, Completely monotonic functions of
positive order and asymptotic expansions of the logarithm of Barnes double
gamma function, \textit{J. Math. Anal. Appl.} \textbf{355} (1) (2009) 33--40.

\bibitem{Lin14} L. Lin, On Stirling's formula remainder,
\textit{Appl. Math. Comput.} \textbf{247} (2014) 494--500.

\bibitem {Magnus66}W. Magnus, F. Oberhettinger, R. P. Soni, \textit{Formulas and
Theorems for the Special Functions of Mathematical Physics}, Springer-Verlag,
Berlin, 1966: 11.

\bibitem {Mortici11}C. Mortici, On the monotonicity and convexity of the remainder
of the Stirling formula, \textit{Appl. Math. Lett.} \textbf{24} (6) (2011) 869--871.

\bibitem {Qi17}F. Qi, B.-N. Guo, Integral representations and complete
monotonicity of remainders of the Binet and Stirling formulas for the gamma
function, \textit{RACSAM} \textbf{111} (2017) 425--434.

\bibitem {Sevli11}H. \c{S}evli, N. Batir, Complete monotonicity results for some
functions involving the gamma and polygamma functions, \textit{Math. Comput.
Modelling} \textbf{53} (9-10) (2011) 1771--1775.

\bibitem {Shi06}X. Shi, F. Liu, M. Hu, A new asymptotic series for the gamma
function, \textit{J. Comput. Appl. Math.} \textbf{195} (1) (2006) 134--154.

\bibitem {Widder46}D. V. Widder, \textit{The Laplace Transform}, Princeton University
Press, Princeton, 1946.


\end{thebibliography}
\end{document}